\documentclass[english]{amsart}

\usepackage{babel}
\usepackage{amstext}
\usepackage{amsmath}
\usepackage{amsfonts}
\usepackage{latexsym}
\usepackage{ifthen}
\usepackage{xypic}
\xyoption{all}
\newcommand{\sO}{{\mathcal O}}

\newcommand{\sE}{{\mathcal E}}

\newcommand{\PN}{{\mathbb P}}
\newcommand{\PF}{{\mathbb F}}
\newcommand{\KQ}{{\mathbb Q}}
\newcommand{\KZ}{{\mathbb Z}}
\newcommand{\KR}{{\mathbb R}}
\newcommand{\Pic}{{\rm Pic}}
\newcommand{\lra}{\longrightarrow}

\newcommand{\Bl}{{\rm Bl}}

\newcommand{\X}{{\mathcal X}}
\newcommand{\Y}{{\mathcal Y}}
\newcommand{\C}{{\mathcal C}}

\newtheorem{lemma1}[equation]{}

\newenvironment{lemma}{\begin{lemma1}{\bf Lemma.}}{\end{lemma1}}

\newenvironment{example}{\begin{lemma1}{\bf Example.}\rm}{\end{lemma1}}
\newenvironment{theorem}{\begin{lemma1}{\bf Theorem.}}{\end{lemma1}}

\newenvironment{Mconjecture}{\begin{lemma1}{\bf
      Generalized Mukai Conjecture.}}{\end{lemma1}}
\newenvironment{theorem2}[1]{\begin{lemma1}{\bf Theorem [#1].}}{\end{lemma1}}
\newenvironment{proposition}{\begin{lemma1}{\bf Proposition.}}{\end{lemma1}}
\newenvironment{corollary}{\begin{lemma1}{\bf Corollary.}}{\end{lemma1}}
 
\newenvironment{remark}{\begin{lemma1}{\bf Remark.}\rm}{\end{lemma1}}

\begin{document}

\title {On the Picard number of almost Fano threefolds with pseudo-index $>1$}
\author[C. Casagrande]{Cinzia Casagrande}
\address{C. Casagrande -
Universit\`a di Pisa - Dipartimento di Matematica - Largo B.\
Pontecorvo, 5 -  I-56127 Pisa, Italy} 
\email{casagrande@dm.unipi.it}
\author[P. Jahnke]{Priska Jahnke} 
\address{P. Jahnke - 
Mathematisches Institut - Universit\"at Bayreuth - D-95440
  Bayreuth, Germany} 
\email{priska.jahnke@uni-bayreuth.de}
\author[I. Radloff]{Ivo Radloff}
\address{I. Radloff - 
Mathematisches Institut - Universit\"at Bayreuth - D-95440
  Bayreuth, Germany} 
\email{ivo.radloff@uni-bayreuth.de}
\thanks{The two last named authors gratefully acknowledge support by the Schwerpunkt program {\em Globale Methoden in der komplexen Geometrie} of the Deutsche Forschungsgemeinschaft.}
\date{1st March 2006}
\subjclass[2000]{14J30, 14J45, 14E30, 14M25}
\maketitle

\section{Introduction}
\noindent A complex projective variety $X$ 
 is called {\em almost Fano}
if it is normal,
Gorenstein, and its anticanonical bundle $-K_X$
is big and nef. 
We assume moreover that $X$ has at most canonical singularities. 

When $-K_X$ is ample, $X$ is a Fano variety. In general, 
by the base point free theorem, some
multiple $|{-}mK_X|$ is base point free, defining a birational
morphism  \[\psi\colon X \lra Y\]
to some Gorenstein Fano variety, again with at most canonical
singularities. Curves contracted by $\psi$ are exactly curves of
anticanonical degree zero. We say that $Y$ is
an {\em anticanonical model} of $X$. Both $X$ and $Y$ are rationally
connected (\cite{KMM}, \cite{zhang}).

An important numerical invariant of $X$ is its \emph{Picard number} $\rho(X)$;
one has $\Pic(X)\cong
H^2(X,\mathbb{Z})$, due to
the Kawamata-Viehweg vanishing
theorem together with the exponential sequence.
Hence  $\rho(X)$ 
coincides with the second Betti number $b_2(X)$. 

In the study of Fano varieties, a relevant role is played by the
index, which is the divisibility of the anticanonical bundle in the
Picard group. Similarly, if $X$ is an almost Fano variety,
one can define its \emph{index} $r_X$ as the divisibility of $-K_X$ in
$\Pic(X)$, and  
its
\emph{pseudo-index} $\iota_X$ as
\[\iota_X := \min\{d > 0 \mid \exists \mbox{ rational curve } C:
-K_X.C = d\}.\]  
Notice that $r_X\mid \iota_X$. As in the Fano case, one expects that
almost Fano varieties with large index or pseudo-index are simpler.

\vspace{0.2cm}

In dimension three, smooth Fano threefolds are classified. Of course
the class of almost Fano threefolds is much larger. In \cite{JPR}
smooth almost Fano threefolds with Picard number two and such that
$\psi$ is divisorial are classified.

If $Y$ is a Gorenstein Fano threefold with at most canonical
singularities,  by results of Kawamata and Reid (\cite{Kaw},
\cite{Reid}) there exists a {\em partial crepant resolution}
 \[\psi\colon X \lra Y,\]
where $X$ is an almost Fano threefold with at most
terminal $\KQ$--factorial singularities, and $Y$ is an anticanonical
model of $X$. Hence singular Fano threefolds and almost Fano threefolds
(with mild singularities) are closely related.

\vspace{0.2cm}

In the present paper we 
study almost Fano threefolds $X$ with at most canonical singularities
and with $\iota_X>1$, namely we assume that $X$ contains no rational
curves of anticanonical degree one.

We first give a birational description of such $X$ in Proposition
\ref{structure}, under the additional assumption that $X$ is terminal and 
$\mathbb{Q}$-factorial. This description, together with
some results 
by Shin and by Chen and Tseng, allows us to characterize
almost Fano threefolds
  with at most canonical singularities and $\iota_X \geq 3$: they
 are just $\PN_3$, quadrics, and resolutions of quadrics (see
 Proposition \ref{description}).

\vspace{0.2cm}

\noindent {\sc Bounding the Picard number.}
It is well known that there exist only finitely many deformation
families of smooth Fano varieties for fixed dimension $n$, 
and the same holds
true in the canonical case at least in dimension three
(\cite{Borisov}, \cite{McK} and \cite{KMMT} for the almost Fano
case). We may therefore ask for (effective) bounds of 
their numerical invariants. 

After the classification, we know that a smooth Fano threefold $Y$
 has $\rho(Y) \le 10$. 
The only case with Picard number $10$ is the product $S
\times \PN_1$, $S$ a Del Pezzo surface of degree
$1$. Bounds for the Picard number of singular Fano
threefolds, or almost Fano threefolds, are still unknown. 
In higher dimensions,
 the maximal Picard number of a smooth Fano variety is
also unknown (even in dimension 4).

In our situation we obtain:

\begin{theorem} \label{almFano}
 Let $X$ be an almost Fano threefold
  with at most canonical singularities and $\iota_X > 1$. 

Then $\rho(X) \le 10$, and equality
holds if and only if $X$ is smooth and there exists a finite sequence of flops $\,X
\dasharrow \Bl_{p_1, \dots, p_8}(\PN(\sO_{\PN_2} \oplus
  \sO_{\PN_2}(3)))$.

If $X$ is Fano, then $\rho(X) \le 3$, and equality holds if and only if $X \simeq
\PN_1 \times \PN_1 \times \PN_1$.
\end{theorem}

\vspace{0.2cm}

In the Fano case, there is a conjectural
relation between the Picard number and the pseudo-index:

\begin{Mconjecture} 
Let $Y$ be a smooth Fano variety. Then
 \[\rho(Y)(\iota_Y-1) \le \dim(Y),\quad\text{with equality if and only
   if}\quad Y\simeq (\PN_{\iota_Y-1})^{\rho(Y)}.\]
\end{Mconjecture}
This was conjectured by Mukai \cite{Mu}
in a weaker form, and then studied in \cite{WisM}, \cite{BCDD},
\cite{occhettaGM}, and  \cite{Ca}. Up to now, the conjecture has been
proven
for smooth Fano varieties of dimension at most 5, and for Gorenstein
and  $\KQ$--factorial
toric Fano varieties of arbitrary
dimension. 
Observe that the conjecture is trivial if $\iota_Y = 1$.

Applying our previous results, we obtain:

\begin{theorem} 
\label{Mukcan}
Let $Y$ be a Gorenstein Fano threefold
  with at most canonical singularities. Then the generalized Mukai
  conjecture holds, i.e.,
 \[\rho(Y)(\iota_Y-1) \le \dim(Y)\quad\text{with equality if and only
   if}\quad Y\simeq (\PN_{\iota_Y-1})^{\rho(Y)}.\]
\end{theorem}

\vspace{0.2cm}

\noindent {\sc The toric case.}
In any 
dimension $n$, Gorenstein toric Fano varieties and almost Fano
toric varieties  are in a finite number, and
always have at most
canonical singularities.

More precisely, Gorenstein
toric Fano varieties of dimension $n$ are in bijection with a special
class of polytopes in $\KR^n$, called \emph{reflexive} polytopes, see
\cite{batdual}. Reflexive polytopes have integral vertices
(i.e.\ in $\KZ^n$), and the origin is their unique interior integral
point. 

If $P$ is a reflexive polytope, we denote by $Y_P$
the corresponding Gorenstein Fano
variety (the fan of $Y_P$ is given by the cones over
the faces of $P$). Then we have:
\begin{enumerate}
\item 
the rank of the divisor class group of $Y_P$ is equal to 
the number of vertices of $P$ minus $n$;
\item the maximal Picard number of toric almost Fano  varieties
whose anticanonical
  model is $Y_P$ is $|P\cap \KZ^n|-n-1$.
\end{enumerate}

Reflexive polytopes have been classified up to dimension four by Kreuzer
and Skarke \cite{PALP}, by means of a computer program. They are 4319
in dimension three, and almost half a billion in dimension four. The
maximal number of integral points are respectively 39 and 680. By 2.)
above, this implies that
if $X$ is an almost Fano toric variety of dimension $n$,
  then 
$$\rho(X)\le\begin{cases} 35 \quad &\text{ if } n=3,\\
                          675 \quad &\text{ if } n=4. 
\end{cases} $$

\vspace{0.2cm}

In arbitrary dimension, the maximal number of integral points of a
reflexive polytope is not known, and sharp bounds on the Picard number
are known only under some additional condition on the singularities.

More precisely, 
if $Y$ is a $\KQ$-factorial
Gorenstein toric Fano variety, it is shown in \cite{Ca} that
$\rho(Y)\le 2n$ if $n$ is even,  $\rho(Y)\le 2n-1$ if $n$ is
  odd.

As seen in 1.),
in the Fano case these combinatorial techniques 
allow to compute the rank of the divisor class group,
so in the non $\KQ$-factorial
case they never give a sharp bound on the Picard number.
In particular, 
Theorem \ref{Mukcan} is new also
in the toric case, at least to our knowledge.

For the almost Fano case, it is shown in \cite{nill}, Corollary 6.3
that
if $X$ is a toric almost Fano of dimension $n$ whose anticanonical
model has at most terminal singularities, then $\rho(X)\le
2^{n+1}-n-2$, and that bound is sharp.

For toric varieties, the bound of Theorem \ref{almFano} can be sharpened:

\begin{proposition} 
\label{toric}
Let $X$ be a toric almost Fano
  threefold  with $\iota_X > 1$. 

Then $\rho(X) \le 8$, and equality can only happen if $X$ is smooth and there
exists a finite sequence of flops  $X\dasharrow \Bl_{p_1, \dotsc,
  p_6}(\PN(\sO_{\PN_2} \oplus 
   \sO_{\PN_2}(3)))$.
\end{proposition}

\section{Almost Fano Threefolds with pseudoindex $>1$}
\setcounter{equation}{0}

\noindent Let $X$ be a Gorenstein almost Fano threefold.
The highest power $({-}K_X)^3 > 0$ is always
even and called {\em anticanonical degree } of $X$. 
Assume that $X$ has canonical singularities, and let
 \[\psi\colon X \lra Y\]
be an anticanonical model of $X$. Since
$\psi$ is crepant, we have
 \[({-}K_X)^3 = ({-}K_Y)^3.\]
Moreover $r_Y=r_X$ and $\iota_Y=\iota_X$.

Since $-K_X=r_X L$, 
the intersection number of $-K_X$ with any rational curve is
divisible by $r_X$, hence $r_X\mid \iota_X$, in particular
 \[r_X \le \iota_X.\] 
As soon as a line exists in $X$, Fano index and pseudo-index coincide
(by a {\em line} we denote a rational curve $C$, such that $L.C = 1$,
cf.\ \cite{Kollar}). So for example on a smooth Fano threefold 
$X$ with
$\rho(X) = 1$ we 
always have $\iota_X = r_X$, since the existence of lines is known. In
general, 
both notions do not necessarily coincide (cf. Lemma~\ref{index}).

\vspace{0.2cm}

Since $K_X$ is not nef, there always exists an \emph{elementary
extremal contraction}
 \[\phi\colon X \lra Z,\]
namely $\phi$ is surjective with connected fibers, $Z$ is normal,
$\rho(X)-\rho(Z)=1$, and $-K_X$ is $\phi$-ample. 
Such
contractions are classified by Mori in the smooth case (\cite{Mori}),
and Cutkosky for terminal and $\KQ$-factorial
singularities (\cite{Cu}). We call an
extremal contraction $\phi$ to be of {\em fiber type}, if $\dim Z <
\dim X$ (including the case $\dim Z = 0$).

\vspace{0.2cm}

Let $X$ be a Gorenstein almost Fano threefold, with at most terminal
singularities. By a \emph{smoothing} of $X$ we mean a flat projective
morphism 
$$\pi\colon\X\lra\Delta$$
onto the unit disc $\Delta$, with $\X$ a reduced and irreducible
complex space, such that $\X_0\simeq X$, while $\X_t$ is 
a smooth almost Fano threefold for $t\neq 0$. 
By \cite{Nam} and \cite{Mi}
a  smoothing always
exists when $X$ is either $\mathbb{Q}$-factorial or Fano. 
In this last
case, $\X_t$ is Fano too. The numerical invariants of $X$ and $\X_t$
are related as follows.
\begin{theorem2}{\cite{smoothing}} \label{Smoothing}
Let $X$ be a Gorenstein almost Fano threefold with at most terminal
singularities, and let $\pi\colon\X\to\Delta$ be a smoothing of
$X$. Then: 
$$(-K_{\X_t})^3=(-K_X)^3,\quad \rho(\X_t)=\rho(X),\quad
 r_{\X_t}=r_X\leq\iota_X\leq\iota_{\X_t}.$$
\end{theorem2}
\begin{proof}
It is shown in \cite{smoothing} that $\X$ is Gorenstein, and that
there is an isomorphism $\Pic(\X_t)\simeq\Pic(X)$ preserving the
canonical class. This implies $(-K_{\X_t})^3=(-K_X)^3$,
$\rho(\X_t)=\rho(X)$, and $r_{\X_t}=r_X$.

 Let $\C_t$ be a rational curve in the
general $\X_t$, $t \not= 0$, such that $-K_{\X_t}.\C_t =
\iota_{\X_t}$. Degenerate $\C_t$ and let $\C_0 \subset X$ be the limit
curve. Then $-K_X.\C_0 = \iota_{\X_t}$ as well, but $\C_0$ might be reducible. This shows
 $\iota_X \le \iota_{\X_t}$.
\end{proof}

\vspace{0.2cm}

We need the following result by Prokhorov:

\begin{theorem2}{\cite{Prokh}} \label{bound}
Let $Y$ be a Gorenstein Fano threefold with at most canonical
singularities. Then $({-}K_{Y})^3 \le 72$ and equality
holds if and only if $Y$ is one of the weighted projective spaces $\PN(1^3,3)$ or
$\PN(1^2,4,6)$. \end{theorem2}

\begin{remark} \label{rembound}
\begin{enumerate}
 \item If $Y$ has only terminal singularities, then $({-}K_Y)^3 \le
   64$. This follows from the existence of a smoothing (\cite{Nam}) in
   this case.
 \item The same bound as in \ref{bound} holds for Gorenstein almost Fano threefolds with
   at most canonical 
singularities. 
 \item A crepant resolution of $\PN(1^3,3)$ is the smooth almost Fano
   threefold $\PN(\sO_{\PN_2} \oplus \sO_{\PN_2}(3))$, which has
   pseudo-index $2$. The pseudo-index of $\PN(1^2,4,6)$ is $1$.
\end{enumerate}
\end{remark}

\begin{example} \label{surfaces}
Smooth almost Fano surfaces can easily be classified, they are
 \begin{enumerate}
   \item $\PN_2$ blown up in $r \le 8$ points,
   \item $\PN_1 \times \PN_1$,
   \item the second Hirzebruch surface $\PF_2$.
 \end{enumerate}
In particular, if $S$ is minimal, then $S \simeq \PN_2$, $\PN_1 \times
\PN_1$ or $\PF_2$. Note that there are conditions on the position of
the points: for example we may allow three points on a line in $\PN_2$, but not
$4$ to ensure that $-K_S$ remains nef.

We may of course blow up $\PF_2$ in a point $p$ not lying on the minimal
section to obtain another almost Fano surface $S = \Bl_p(\PF_2)$. Let $f
\subset \PF_2$ be the fiber containing $p$. Then the strict transform
of $f$ in $S$ becomes a $(-1)$--curve. Contracting this curve to a
point $p'$, we obtain another Hirzebruch surface, namely $\PF_1 \simeq
\Bl_q(\PN_2)$. Hence $S \simeq \Bl_{p',q}(\PN_2)$. This shows that the
above list is indeed complete.
\end{example}

Before we come to the proof of Theorem~\ref{almFano}, we need some
general results on possible elementary
extremal contractions of an almost Fano
threefold with pseudo-index $\iota_X>1$, under some additional
assumptions on the singularities of $X$ (for similar results in the
smooth case compare \cite{DPS}).

\begin{lemma} \label{conicbdl}
Let $X$ be an
almost Fano threefold with at most terminal and $\KQ$-factorial
singularities and $\iota_X > 1$.
Suppose that $X$ admits an elementary extremal
contraction
 \[\phi\colon X \lra S\]
onto a surface $S$. Then
$S$ is a
smooth surface with $-K_S$ big and nef, and  $X = \PN(\sE)$ for some
rank $2$ vector bundle 
$\sE$ on $S$.
\end{lemma}

\begin{proof}
By classification, $S$ is smooth and $\phi$ is a conic
bundle  (see \cite{Cu}, Theorem 7). 
Since $\iota_X > 1$ by assumption, $\phi$ has no singular fibers.
Then $X$ is smooth and it remains to show $-K_S$ big and
nef. By \cite{DPS}, Proposition~3.1, $-K_S$ is nef and $X = \PN(\sE)$ for
some rank two vector bundle $\sE$ on $S$. By \cite{AG5}, Proposition 7.1.8, we have
 \[-4K_S = \phi_*({-}K_X)^2 + \Delta,\]
where $\Delta$ is the discriminant of the conic bundle, hence $\Delta
= 0$. By the Riemann--Roch theorem and Kawamata--Viehweg vanishing,
$|{-}K_X|$ is non--empty and the sections cover $X$. Hence $-4K_S = \phi_*C$ for some complete intersection curve $C
= H_1 \cap H_2$ with $H_i \in |{-}K_X|$ is effective and
irreducible. Moreover, moving for example $H_2$ in $X$, we find
$(\phi_*C)^2 > 0$. This implies $-K_S$ is also big as claimed. 
\end{proof}

\begin{lemma} \label{notmin}
Let $X$ be as in Lemma \ref{conicbdl}. Assume that the surface $S$ is 
not minimal. Then there exists a flop diagram
 \begin{equation} \label{flop}
   \xymatrix{X \ar@{-->}[rr]^{\chi} \ar[dr] & & X^+ \ar[dl]\\
             & Y &}
 \end{equation}
over the anticanonical model $Y$ of $X$, such that $X^+$ admits a
birational elementary extremal
contraction $\phi^+\colon X^+ \to Z^+$. Moreover, $X^+$ is again an
almost Fano threefold with at most terminal and $\KQ$-factorial
singularities, and
 \[\iota_{X^+} = \iota_X, \quad \rho(X^+) = \rho(X) \quad \mbox{ and }
 \quad ({-}K_{X^+})^3 = ({-}K_X)^3.\] 
\end{lemma}

\begin{proof}
Since $S$ is not minimal by assumption, there exists a $(-1)$--curve
$C \subset S$, i.e., 
 \[C \simeq \PN_1, \quad C^2 = -1 \quad \mbox{ and } -K_S.C = 1.\] 
Let $X = \PN(\sE)$ for some rank $2$ vector bundle $\sE$ on $S$ and twist by a line
bundle, such that 
 \[\sE|_C = \sO_{\PN_1} \oplus \sO_{\PN_1}(a), \quad \mbox{ for some }
 a \ge 0.\]

Let $\tilde{S} = \phi^{-1}(C) = \PN(\sE|_C)$ and $\tilde{C} \to C$ the section in
$\tilde{S}$ over $C$ corresponding to the projection of $\sE|_C \to
\sO_{\PN_1} \to 0$. Let $\sO(1)$ be the tautological line bundle on
$X$ as usual. Then 
 \[-K_X.\tilde{C} = \sO(2).\tilde{C} -K_S.C - (\det \sE).C = 1-a\geq 0,\]
so $a=0,1$. By assumption $\iota_X>1$,
hence $a = 1$ and $\tilde{C}$ is an anticanonically trivial curve in
$X$. Moreover, 
 \[\tilde{S} \simeq \PF_1 = \PN(\sO_{\PN_1} \oplus
 \sO_{\PN_1}(1)),\]
where $-K_X|_{\tilde{S}} = \sO(2)$ is $2$ times the corresponding
tautological line bundle of $\PF_1$ and $\tilde{C}$ is the minimal
section. We have $\sO_X(\tilde{S})|_{\tilde{C}} = \sO_{\PN_1}(-1)$,
since $C$ is a $(-1)$--curve by definition, and
$N_{\tilde{C}/\tilde{S}} = \sO_{\PN_1}(-1)$, since $\tilde{C}$ is the
minimal section of $\tilde{S}$. This implies the splitting of the
normal bundle sequence
 \[0 \lra N_{\tilde{C}/\tilde{S}} \lra N_{\tilde{C}/X} \lra
 N_{\tilde{S}/X}|_{\tilde{C}} \lra 0\]
and therefore the splitting type of $N_{\tilde{C}/X}$ is
$(-1,-1)$. Blowing up the curve $\tilde{C}$ in $X$, the exceptional
divisor will be $\PN_1 \times \PN_1$, and we may blow down in the
other direction onto some variety $X^+$
 \[\xymatrix{& \Bl_{\tilde{C}}(X) \ar[dl]_{\pi} \ar[dr]^{\pi^+} &\\
             X \ar@{-->}[rr]^{\chi} && X^+}\] 
where the rational map $\chi$ is the flop (\ref{flop}), since
$-K_X.\tilde{C} = 0$. In particular, $X^+$ is again an almost Fano
threefold with at most terminal 
and $\KQ$-factorial singularities (cf.\ \cite{KoMo}, \S 6.2). 

It remains to show $\iota_{X^+} > 1$ and that $X^+$ now admits a
birational elementary
extremal contraction. First note that $\chi$ maps the surface
$\tilde{S}$ onto some 
$S^+ \simeq \PN_2$ in $X^+$ with normal bundle $N_{S^+/X^+} =
\sO_{\PN_2}(-1)$. This means we may contract $S^+$ to a smooth point,
which gives the map $\phi^+$ in the lemma.

In order to see $\iota_{X^+} = \iota_X$ let $C^+ \subset X^+$ be any
rational curve. Let $C_0 \subset X$ be the strict transform of $C^+$,
i.e., $C_0 = 
\pi(C_0^+)$ for some section $C_0^+ \subset \Bl_{\tilde{C}}(X)$ over
$C^+$. Then
 \[-K_X.C_0 = \pi^*(-K_X).C_0^+ = (\pi^+)^*(-K_{X^+}).C_0^+ = -K_{X^+}.C^+.\]
\end{proof}

\begin{lemma} \label{birat}
Let $X$ be an
almost Fano threefold with at most terminal and $\KQ$-factorial
singularities and $\iota_X > 1$.
Suppose that $X$ admits a birational elementary extremal
contraction
$$\phi\colon X \lra Z. $$
Then $Z$ is an almost Fano threefold with at most
terminal and $\KQ$-factorial singularities, $\iota_Z>1$,
 and $\phi$ is the blowup of a smooth point.
\end{lemma}

\begin{proof}
Since $X$ is Gorenstein with terminal 
and $\KQ$-factorial singularities, Cutkosky's
classification applies (\cite{Cu}): $\phi$ is divisorial and either contracts the
exceptional divisor $E$ onto a local complete intersection curve contained in the
smooth locus of $Z$, or to a point. The first case is impossible,
since then a general fiber of $E$ is a rational curve $f$ with $-K_X.f
= 1$. 

If $E$ is mapped to a point, then the pair $(E, E|_E)$ is one of
$(\PN_2, \sO(-1))$, $(\PN_2, \sO(-2))$ or $(Q, \sO(-1))$, where $Q$ is
either a smooth quadric or the quadric cone. The latter two cases are
again impossible, since they provide curves of anticanonical degree
$1$ in $X$. Hence  $(E, E|_E) = (\PN_2, \sO(-1))$ and in particular $E$
is mapped to a smooth point of $Z$ and $Z$ is again Gorenstein with at
most terminal and $\KQ$-factorial singularities.

Let $C \subset Z$ be any rational curve and $\tilde{C} \to C$ a
section in $X$. If $C$ does not contain the point $p = \phi(E)$, then
$-K_Z.C = -K_X.\tilde{C}$. If $p \in C$, then 
 \begin{equation} \label{iota2}
   -K_Z.C = -K_X.\tilde{C} + 2E.\tilde{C} \geq -K_X.\tilde{C}+2,
 \end{equation}
since $\tilde{C} \not\subset E$. This shows $\iota_Z >1$ 
and that $-K_Z$ is
nef. Finally $({-}K_Z)^3 = ({-}K_X)^3 + 8 > 0$, hence $-K_Z$ is
also big. 
\end{proof}

\begin{proposition}
\label{structure}
Let $X$ be an
almost Fano threefold with at most terminal and $\KQ$-factorial
singularities and $\iota_X > 1$. Then there exists a sequence
\begin{equation} \label{fact}
  \xymatrix{X\ar@{-->}[r]^{\xi}& X_0 \ar[r]^{\phi_0} & X_1 \ar[r]^{\phi_1} &
   \cdots \ar[r]^{\phi_{m-1}}  &  X_m \ar[r]^{\phi} & W}
  \end{equation}
where:
\begin{enumerate}
\item each $X_i$ is an
almost Fano threefold with at most terminal and $\KQ$-factorial
singularities and $\iota_{X_i}>1$;
\item $\xi$ is a finite sequence of flops and $\iota_X=\iota_{X_0}$;
\item each $\phi_i$ is a blowup in a smooth point;
\item $\phi$ is an elementary extremal 
contraction of fiber type;
\item $W$ is one of the following: 
a point, $\PN_1$, $\PN_2$, $\PN_1\times\PN_1$, or $\mathbb{F}_2$;
\item if $W$ is $\PN_1\times\PN_1$ or $\mathbb{F}_2$, then $m=0$.
\end{enumerate}
\end{proposition}
\begin{remark}
When $\dim W=2$, then $X_m$ is a $\PN_1$-bundle over $W$ by Lemma
\ref{conicbdl}, so $X$ is smooth.

When $W\simeq \PN_1$, then $X_m\to W$ is a del Pezzo fibration.
 If $F$ is a general fiber, we have 
$-K_{X_m}|_F = -K_F$, so $\iota_F > 1$. This means that
$F$ is either $\PN_2$ or $\PN_1\times \PN_1$.
\end{remark} 
\begin{proof}[Proof of Proposition \ref{structure}]
Let $X$ be
 an
almost Fano threefold with at most terminal and $\KQ$-factorial
singularities and $\iota_X > 1$.
As seen previously, there exists an elementary extremal contraction
 \[\phi\colon X \lra Z.\]
By \cite{Cu}, we have the following possibilities:
\begin{enumerate}
 \item $\dim Z = 0$. 
 \item $\dim Z = 1$, hence $Z\simeq\PN_1$. 
 \item $\dim Z = 2$ and $\phi$ is a conic bundle. In this case $\phi$
   is a $\PN_1$-bundle over a smooth almost Fano surface $Z$ by
   Lemma~\ref{conicbdl} and we are in one of the following two cases:
   \begin{enumerate}
    \item $Z$ is minimal, hence $Z$ is one of the surfaces
$\PN_2$, $\PN_1\times\PN_1$, or $\PF_2$ (cf. example \ref{surfaces}).
     \item $Z$ is not minimal. In that case we apply Lemma~\ref{notmin}
      performing a flop $X \dasharrow X^+$, where $X^+$
      is another almost Fano threefold with at most terminal
      and $\KQ$-factorial
singularities, such that there exists an elementary extremal contraction
      $\phi^+\colon X^+ \to Z^+$ with $\dim Z^+ = 3$. We have $\iota_X
      = \iota_{X^+}$, $\rho(X) = \rho(X^+)$ and $({-}K_X)^3 =
      ({-}K_{X^+})^3$. 
   \end{enumerate}
 \item $\dim Z = 3$. Then by Lemma~\ref{birat}, $\phi$ is the blowup
   of a smooth point, $Z$ is again almost Fano with at most
   terminal and $\KQ$-factorial singularities and $\iota_Z > 1$.
\end{enumerate}
This means that we are done if either $\dim Z \le 1$ or $\dim Z = 2$ and
$Z$ is minimal. As long as this is not the case, we proceed as
follows: first blow down divisors to points as long as possible. If we
then end up in case 1.), 2.) or 3.) (i), we stop. If we end up in case 3.) (ii), we perform a flop and
start again blowing down. This process is finite, since on the
one hand the anticanonical degree remains stable in case of a flop and
increases by $8$ if we blow down a divisor to a smooth point. On the
other hand, 
$(-K_X)^3 \le 72$ as seen in Theorem~\ref{bound}. 

\vspace{0.2cm}

Hence we finally end up with the following picture:
\begin{equation}\label{unordered}
  \xymatrix{X = {X}'_0 
\ar@{-->}[r]^{{\phi}'_0} & {X}'_1
\ar@{-->}[r]^{{\phi}'_1}  &
   \cdots \ar@{-->}[r]^{\phi'_{n-1}}  &  X'_n \ar[r]^{\phi} & W}
\end{equation}
where $\phi$ is an elementary extremal 
contraction of fiber type, $W$ is as in 5.), and each
$\phi'_i$ is either a flop or a blowup of a smooth point.
By Lemma~\ref{birat} and
Lemma~\ref{notmin}, each $X'_i$ is again an almost Fano threefold with
at most terminal and $\KQ$-factorial singularities and $\iota_{X'_i} > 1$. 

\vspace{0.2cm}

Suppose that for some index $i\in\{0,\dotsc,n-2\}$ we have 
 $$\xymatrix{{X'_{i}} \ar[r]^{\phi'_{i}} & {X'_{i+1}} 
\ar@{-->}[r]^{\phi'_{i+1}} & {X'_{i+2}}}$$
with $\phi'_{i}$ a blowup of a smooth point $p\in X'_{i+1}$,
and $\phi'_{i+1}$ a flop. 
Since $-K_{X'_{i}}$ is nef, 
the point $p$ can not lie on any
anticanonically trivial curve, in particular it will not lie on the
exceptional locus of the flop. So we can first blowup the image of $p$
in $X'_{i+2}$ and then perform the flop; in this way 
 we get a new factorization of 
$\phi'_{i+1}\circ\phi'_{i}$ as
$$\xymatrix{{X'_{i}} \ar@{-->}[r]^{\phi''_{i}} & {X''_{i+1}} 
\ar[r]^{\phi''_{i+1}} & {X'_{i+2}}}$$
where by
Lemma~\ref{notmin}, $X''_{i+1}$ is again 
an almost Fano threefold with
at most terminal and $\KQ$-factorial
singularities and $\iota_{X''_{i+1}} > 1$.

Iterating this procedure and  renaming, we are reduced to a sequence
$$
  \xymatrix{X\ar@{-->}[r]^{\xi}& X_0 \ar[r]^{\phi_0} & X_1 \ar[r]^{\phi_1} &
   \cdots \ar[r]^{\phi_{m-1}}  &  X_m \ar[r]^{\phi} & W}
$$
where $\xi$ is a sequence of flops, and 
each $\phi_i$ is a
blowup of a smooth point $p_{i+1}\in X_{i+1}$. 

\vspace{0.2cm}

Suppose that $m>0$ and $W=\mathbb{F}_2$ or $W=\PN_1\times\PN_1$.
We want to see that these cases can be
reduced to the case $W=\PN_2$ by a sequence of flops. 

Assume $W =\PF_2$. We have $X_{m-1} =
\Bl_p(X_m)$ for some 
point $p$. The strict transform of the fiber containing $p$ of the
$\PN_1$--bundle $X_m \to W$ becomes anticanonically trivial in
$X_{m-1}$ with normal bundle of type $(-1,-1)$. Flopping that curve we
obtain a $\PN_1$--bundle $X'_{m-1}$ over $\Bl_p(\PF_2)$ (we call the image of
$p$ in $W$ again $p$; this is the reversed construction of the flop
over a $({-}1)$--curve in Lemma~\ref{notmin}). 

But $\Bl_p(\PF_2)$ now admits a second $({-}1)$--curve, namely the strict transform of the
fiber of $\PF_2$ containing $p$. Blowing this curve down, we obtain
$\PF_1$ (compare Example~\ref{surfaces}). On the other hand, this second 
$({-}1)$--curve gives rise to another flop of $X'_{m-1}$, we call the
resulting threefold $X''_{m-1}$. Blowing down the exceptional divisor
in $X''_{m-1}$ provided by the flop $X'_{m-1} \dasharrow X''_{m-1}$ we obtain a $\PN_1$--bundle $X''_m$ over
$\PF_1$. Flopping over the minimal section of $\PF_1$, which is now a
$({-}1)$--curve, we finally end up with a $\PN_1$--bundle $X''_{m+1}$
over $\PN_2$. In a diagram:

\[\xymatrix{& X''_{m-1} \ar[r]^{\Bl} & X''_m \ar@{-->}[r]^{\text{\em flop}}
  \ar'[d][dd] &
  X'''_m \ar[r]^{\Bl} & X'''_{m+1} \ar[dd]\\ 
            X'_{m-1} \ar@{<-->}[rr]^{flop} \ar@{<-->}[ur]^{
\text{\em flop}}
            \ar[dd] && X_{m-1} \ar[r]^{\Bl} & X_m \ar[dd]&\\
      &&{\PF_1} \ar'[r][rr] && \PN_2\\
            \Bl_p(\PF_2) \ar[urr] \ar[rrr] &&& \PF_2&}\]

Thus we are reduced to a new sequence as \eqref{unordered}, ending
with $\PN_2$ instead of $\mathbb{F}_2$. Now we repeat the procedure of
ordering the flops and the blowups, and get
the statement. The
same argument applies for $W = \PN_1 \times \PN_1$.
\end{proof}

\begin{proof}[Proof of Theorem~\ref{almFano}, almost Fano case.]
Let $X$ be a Gorenstein almost Fano threefold with at most canonical
singularities and $\iota_X > 
1$. First of all, we reduce to the case where $X$ has at most 
$\KQ$-factorial and terminal singularities. In
fact, by \cite{KoMo}, Theorems 6.23 and 6.25, 
there exist two birational morphisms
$$ X''\stackrel{g}{\lra}X'\stackrel{f}{\lra}X$$
such that $f$ is crepant, $X'$ has at most terminal singularities, $g$ 
is an isomorphism in codimension
1, and $X''$ has at most terminal and
$\KQ$-factorial singularities.

We have $K_{X'}=f^*(K_X)$ and   $K_{X''}=g^*(K_{X'})$, so both $X'$
and $X''$ are almost Fano and
$\iota_{X''}=\iota_{X'}=\iota_X>1$. Moreover
$\rho(X)\leq\rho(X')\leq\rho(X'')$, and $\rho(X)=\rho(X'')$ if and only if
$X\simeq X''$ is already terminal and $\KQ$-factorial.

So we assume that $X$  has at most 
$\KQ$-factorial and terminal singularities.
Applying Proposition \ref{structure}, we get a sequence as
\eqref{fact}, so that
\begin{gather*}
\rho(X) = \rho(X_0)=\rho(X_m)+m=
\rho(W) + (m+1),\\
({-}K_X)^3 =({-}K_{X_0})^3 = ({-}K_{X_m})^3- 8m.\end{gather*}
Observe that $\rho(W)\leq 2$, so if $m=0$ we get $\rho(X)\leq 3$ and
we are done.

We now assume $m>0$ and consider the possible cases for $W$ separately: 

\vspace{0.2cm}

1.) If $W$ is a point, then $X_m$ is a Fano threefold with 
terminal
Gorenstein singularities, hence $({-}K_{X_m})^3 \le 64$ by
Remark~\ref{rembound}. Then   
 \[2 \le ({-}K_X)^3 = ({-}K_{X_m})^3 - 8m \le 64 - 8m,\]
hence $m \le 7$. This gives $\rho(X) \le \rho(W) + 8 = 8$. 

\vspace{0.2cm}

2.) If $W\simeq\PN_1$, then $X_m \to \PN_1$ is a del Pezzo fibration. 
 We claim that $({-}K_{X_m})^3 \le 64$. This implies  $m \le 7$ and 
therefore $\rho(X)\le\rho(W)+8 = 9$.

By \cite{Mi}, there exists a
smoothing $\X \to \Delta$ of $X_m$; 
the general fiber $\X_t$ is a smooth almost Fano with
$(-K_{\X_t})^3=({-}K_{X_m})^3$ and  $\rho(\X_t) = \rho(X_m)=2$
(see Theorem \ref{Smoothing}).

Let $\psi_t$ be the anticanonical map of $\X_t$. If $\psi_t$ is divisorial, then
$({-}K_{X_m})^3 \le 64$ by the classification in \cite{JPR}. If $\psi_t$
is small, then the anticanonical model $\Y_t$ of $\X_t$ is a Fano threefold
with at most terminal singularities. Hence $({-}K_{X_m})^3 =
({-}K_{\Y_t})^3 \le 64$ by Remark~\ref{rembound}. 

\vspace{0.2cm}

3.) If $W = \PN_2$, then by construction $X_m$ is a $\PN_1$--bundle over 
$W$. Since $({-}K_{X_m})^3 \le 72$ by Theorem~\ref{bound}, we find $m\leq
8$. Then
$\rho(X)\leq \rho(W)+9\leq 10$.

\vspace{0.2cm}

Assume now that $\rho(X) = 10$ is maximal for some Gorenstein almost Fano
threefold $X$ with at most canonical singularities and $\iota_X >
1$. Then $X$ has only terminal $\KQ$--factorial
singularities. Consider again the sequence (\ref{fact}) for $X$. 
As seen above, 
it must be $W=\PN_2$, $m=8$ and 
$({-}K_{X_m})^3 =
72$. 
Theorem~\ref{bound} and Remark~\ref{rembound} yield 
$X_m=\PN(\mathcal{O}_{\PN_2}\oplus\mathcal{O}_{\PN_2}(3))$, and the
result follows. 
\end{proof}

\begin{example} \label{Ex10}
We want to construct almost Fano threefolds $X$ with pseudo-index $2$
and Picard number $10$. Let for example $S$ be a del Pezzo surface of degree $1$, i.e., $S = \Bl_{p_1,
  \dots, p_8}(\PN_2)$. Define
 \[X = \PN(\sO_S \oplus \sO_S(-K_S)).\]
Then $-K_X = \sO_X(2)$ is two times the tautological bundle, hence
$r_X = \iota_X = 2$ and $-K_X$ is nef. Moreover $({-}K_X)^3 = 8$ shows $-K_X$ is
also big. Since $\rho(S) = 9$, we have $\rho(X) = 10$ as claimed.

Following the proof of Theorem \ref{almFano}, 
 $X$ should be connected
to a $\PN_1$--bundle over $\PN_2$ by a
sequence of flops and blowups. This can be seen as follows: over each
$(-1)$--curve in $S$ lies an $\PF_1$. Flopping the minimal section,
yields another almost Fano threefold, where we now may contract the
resulting $\PN_2$ to a point. We finally arrive at
 \[\xymatrix{X = X_0 \ar@{-->}[r] & X_8 = \PN(\sO_{\PN_2} \oplus
   \sO_{\PN_2}(3)) \ar[r] & \PN_2 = Z.}\]

Consider on the other hand
 \[X' = \Bl_{p_1, \dots, p_8}(\PN(\sO_{\PN_2} \oplus \sO_{\PN_2}(3))),\]
the blowup of $\PN(\sO_{\PN_2} \oplus \sO_{\PN_2}(3))$ in $8$ general
points. Then $X'$ is another almost Fano 
threefold with pseudo-index $\iota_{X'} = 2$ and maximal Picard number
$\rho(X') = 10$. We may view $X$ and $X'$ as two different crepant
resolutions of the same anticanonical model $Y$.

\vspace{0.2cm}

Note that this is the maximal number of points we may blow up in
either $S$ or the projective bundle $\PN(\sO_{\PN_2} \oplus \sO_{\PN_2}(3))$:
if we blow up one further point, the anticanonical degree decreases
to zero, i.e., that threefold is not almost Fano anymore.
\end{example}

\section{Almost Fano threefolds with high index}
\setcounter{equation}{0}

\noindent We recall two results about Fano threefolds with high index or
pseudo-index by Shin and by Chen and Tseng.

\begin{theorem2}{\cite{Shin}, Theorem 3.9}\label{shin}
Let $Y$ be a Gorenstein Fano threefold with at most canonical
singularities. Then:
\begin{enumerate}
\item
 $r_Y\leq 4$, with equality if and only if
$Y\simeq\PN_3$;
\item $r_Y=3$ if and only if $Y$ is a (possibly
singular) quadric in $\PN_4$.
\end{enumerate}
\end{theorem2}

\begin{theorem2}{\cite{CT}, Corollary 5.2}\label{ct}
Let $Y$ be a Fano threefold with at most canonical
singularities. Then $\iota_Y\leq 4$, with equality if and only if
$Y\simeq\PN_3$.
\end{theorem2}
We give an analogous result about the almost Fano case.

\begin{proposition} 
\label{description}
Let $X$ be an almost Fano threefold with at most canonical
singularities. Then:
\begin{enumerate}
\item $\iota_X\leq 4$, with equality if and only if $X\simeq\PN_3$;
\item
$\iota_X=3$ if and only if $X$ is one of the following:
a (possibly singular)
quadric,
$\PN(\mathcal{O}_{\PN_1}\oplus \mathcal{O}_{\PN_1}(1)^{\oplus 2})$,
or
$\PN(\mathcal{O}_{\PN_1}^{\oplus 2}\oplus \mathcal{O}_{\PN_1}(2))$.
\end{enumerate}
\end{proposition}
\begin{proof}
Let $Y$ be an anticanonical model of $X$, then $\iota_X=\iota_Y$. By
Chen and Tseng's
Theorem~\ref{ct}, we have  $\iota_Y\leq 4$, with equality  if
and only if $Y\simeq\PN_3$. In this case, it must also be
$X\simeq\PN_3$. This gives 1.). 

\medskip

Assume now that $\iota_X=3$. Observe that since $r_X|\iota_X$, we have
$r_X\in\{1,3\}$.

We first show that if $X$ has at most terminal and $\mathbb{Q}$-factorial
singularities, then $r_X=3$. By contradiction, assume $r_X=1$.

There exists a smoothing $\X \to
\Delta$ of $X$ by \cite{Mi}, and
by  Theorem \ref{Smoothing}
the general fiber
$\X_t$ is a smooth almost Fano threefold with $r_{\X_t}=1$,
$\iota_{\X_t}\geq 3$, and $\rho(\X_t)=\rho(X)$.

If $\X_t$ is Fano, then
$\X_t$ is $\PN_3$ or a quadric by \cite{Miy} (or just by
classification), which contradicts $r_{\X_t}=1$. Hence $\X_t$ and $X$
are not Fano.

Let's show that $\rho(\X_t)=2$.
Let $\phi\colon X \to Z$ be any elementary
extremal contraction. By Cutkosky's classification
(see Lemmas \ref{conicbdl} and \ref{birat}), $X$ contains
rational curves of anticanonical degree $1$ or $2$, except if $\phi$ is a
del Pezzo fibration with general fiber $\PN_2$, or $Z$ is a point. 
The latter case can be excluded because $X$ is not Fano,
hence $\rho(\X_t)=\rho(X)=2$.

Let $\psi_t\colon
\X_t \to Y_t$ be an anticanonical model of $\X_t$. Then $\rho(Y_t) =
1$, $r_{Y_t}=1$ and $\iota_{Y_t} \ge 3$. If $\psi_t$ is small, then
$Y_t$ 
is terminal and admits
a smoothing $\Y_{t,s}$, which (again by Theorem 
\ref{Smoothing}) has pseudo-index at least 3 and index 1.
As above, $\Y_{t,s}$  must be either $\PN_3$ or a quadric, and we get
a contradiction.
If $\psi_t$ is divisorial, then from \cite[Table A.2]{JPR} we see
that $Y_t$ is a quadric, which is again
impossible. 

\medskip

Now let $X$  be an almost Fano threefold with at most canonical
singularities and $\iota_X=3$.
As in the proof of Theorem \ref{almFano}, there exists a crepant
birational morphism $f\colon X'\to X$ where $X'$ is an almost Fano
threefold with at most 
terminal and $\mathbb{Q}$-factorial singularities, and
$\iota_{X'}=\iota_X=3$. Then $r_{X'}=3$,
and an
anticanonical model of $X'$ is a quadric by Shin's Theorem \ref{shin}.
Therefore $X'$ is either a quadric, or 
$\PN(\mathcal{O}_{\PN_1}^{\oplus 2}\oplus \mathcal{O}_{\PN_1}(2))$, or
$\PN(\mathcal{O}_{\PN_1}\oplus \mathcal{O}_{\PN_1}(1)^{\oplus 2})$. 
Hence either $X=X'$, or $X$ is a quadric.
\end{proof}
Finally, we give some properties of almost Fano threefolds with $r_X\neq\iota_X$.
\begin{corollary} \label{index}
Let $X$ be an almost Fano threefold with at most terminal and
$\mathbb{Q}$-factorial 
singularities and $r_X\neq\iota_X$. 
Then $\rho(X)=2$, $r_X=1$, $\iota_X=2$, and $X$ is smooth
if and only if $X \simeq \PN_1 \times \PN_2$.
\end{corollary}
\begin{proof}
First note that $X = \PN_1 \times \PN_2$ has Fano--index $r_X = 1$, but
pseudo-index $\iota_X = 2$: the anticanonical divisor is of bidegree
$(2,3)$, hence not divisible in $\Pic(X)$, but $X$ contains no curves
of degree $1$. 

\vspace{0.2cm}

Let conversely $X$ be an almost Fano threefold with at most terminal and
$\mathbb{Q}$-factorial singularities and $r_X\neq \iota_X$.
Then Proposition \ref{description} implies that $\iota_X=2$ and $r_X=1$.

By \cite{Mi} there exists a smoothing $\X \to
\Delta$ of $X$, and 
by Theorem
\ref{Smoothing}
the general fiber
$\X_t$ is a smooth almost Fano threefold with $r_{\X_t}=1$,
$\iota_{\X_t}\geq 2$ and $\rho(\X_t)=\rho(X)$. 
Again we must have $\iota_{\X_t}=2$.

It remains to show $X \simeq \PN_1 \times \PN_2$ for smooth $X$: indeed,
then $2 = \rho(\X_t) = \rho(X)$ completes the proof. So assume $X$
smooth with $\iota_X = 2$ and $r_X = 1$. 
Consider the chain (\ref{fact}) for $X$. We want to see that
 \begin{equation} \label{rxm}
  r_{X_m} =1.
 \end{equation}
Indeed $r_{X_0}=r_X=1$, so if $m=0$ we are done. Assume that $m>0$,
let $i\in\{0,\dotsc,m-1\}$ and let 
$E_i\subset X_i$ be the exceptional divisor of $\phi_i\colon X_i\to X_{i+1}$. Since 
$\Pic(X_i)\simeq\phi_i^*\Pic(X_{i+1})\oplus\mathbb{Z}E_i$ and 
 $-K_{X_i} = \phi_i^*({-}K_{X_{i+1}}) -2E_i$, we have 
 $r_{X_i}=\gcd(2,r_{X_{i+1}})$. On the other hand, any line in $E_i$
 has anticanonical degree 2, so $\iota_{X_i}=2$. 

Therefore for $i<m-1$, $r_{X_i}=1$ implies that $r_{X_{i+1}}$ is odd and at
most 2, i.e. $r_{X_{i+1}}=1$. Thus we get that $r_{X_{m-1}}=1$ and
$r_{X_m}$ is odd. It is then enough to show that $r_{X_m}\neq 3$. 

If $r_{X_m}=3$, then by Proposition
\ref{description} $X_m$ is either a quadric or a resolution of a
quadric. In any case, through every point of $X_m$ there is a smooth
curve of anticanonical degree 3, which implies that $X_{m-1}$ should 
have pseudo-index 1, a contradiction. This completes the proof of (\ref{rxm}).

\medskip

If $X_m$ is Fano, then $1=r_{X_m}<\iota_{X_m}$ implies $X_m
\simeq \PN_1 \times \PN_2$ by a result of Shokurov 
(\cite{Shok}, under the additional assumption that $-K_{X_m}$ is very
ample. But the remaining cases can easily by solved by
classification; compare \cite{AG5}, Theorem 2.4.5 and Theorem
2.1.16 for the respective lists). 
Since the blowup of $\PN_1\times\PN_2$ in a point has pseudo-index 1,
it must be $m=0$ and $X\simeq\PN_1\times\PN_2$. 

We may hence assume
that $X_m$ is not Fano, and show that this gives a contradiction.
Consider the extremal
contraction $\phi\colon X_m\to W$.
Observe that $\rho(W)\in\{1,2\}$ and $\rho(W)=2$ only if $W = \PN_1 \times \PN_1$ or $W =
\PF_2$. We show that in this two cases, the index of $X$ can not be one.

Assume that
$W = \PN_1 \times \PN_1$ and write $X_m = \PN(\sE)$ for some rank $2$ vector bundle $\sE$
on $W$. Twisting by a line bundle we may assume that 
 \[\sE|_{l_1} = \sO_{\PN_1} \oplus  \sO_{\PN_1}(a) \quad \mbox{ and }
 \quad \sE|_{l_2} = \sO_{\PN_1} \oplus  \sO_{\PN_1}(b) \quad \mbox{ for
   some } a,b \ge 0,\] 
where $l_1, l_2 \simeq \PN_1$ are the two rulings of $\PN_1 \times
\PN_1$. Then $\det \sE$ is divisor of bidegree $(b,a)$. Let $C_i$ be
the section in $X$ over $l_i$ corresponding to the projection 
 \[\sE|_{l_i} \lra \sO_{\PN_1} \lra 0.\]
Then $-K_X.C_1 = 2-a$ and $-K_X.C_2 = 2-b$, hence $a,b \in \{0,2\}$
are both even (note that there exists no rational curve $C$ in $X$,
such that $-K_X.C = 1$, and that $-K_X$ is nef). If $\sO(1)$ denotes the tautological line bundle on
$X$, we have $-K_X = \sO(2) \otimes \phi^*\sO(2-b, 2-a)$ is divisible
by $2$, contradicting \eqref{rxm}. If $W = \PF_2$ take a fiber of the ruled
surface $\PF_2$ and the minimal section instead of the rulings. Then
the same argument applies.

\medskip

Hence we have $\rho(W)=1$ and $\rho(X_m)=2$.
Let $\psi\colon X_m\to Y_m$ be an anticanonical model. Since
$\rho(X_m) = 2$ and $\psi$ is non--trivial, we have $\rho(Y_m) = 1$. If $\psi$ is
small, then $Y_m$ is terminal, Fano, with $\iota_{Y_m}>r_{Y_m}$. Then
a smoothing $\Y_{m,t}$ of $Y_m$ exists and is Fano 
with $\rho(\Y_{m,t}) = 1$ and $\iota_{\Y_{m,t}} \ge \iota_{Y_m} >
r_{Y_m} = r_{Y_{m,t}}$. This contradicts \cite{Shok}. Hence $\psi$ is
divisorial.

\medskip

We are left with the following possibilities:
\begin{enumerate}
\item $W\simeq\PN_1$ and the general fiber of $\phi$ is
  $\PN_2$: then from \cite[Table A.2]{JPR} we see that $Y_m$ is a
  quadric hence $r_{X_m}=3$, impossible;
\item $W\simeq\PN_1$ and the general fiber of $\phi$ is
  $\PN_1\times\PN_1$: then by \cite[Table A.2]{JPR} there are four
  possibilities for $X_m$ (N.\ 9, 12, 14, 15). For cases 9, 14, 15 the
  (2) in the column corresponding to $X'$ in the table indicates that
  $\psi$ may be defined by $|{-}\frac{1}{2}K_{X_m}|$, hence $r_{X_m}$
  is divisible by two. For case 12 note that the anticanonical model
  is the cone over the Verones surface, which has index two. So this is impossible. 
\item $W\simeq\PN_2$: then by \cite[Table A.3]{JPR} we see that there
  are four possibilities for $X_m$ (N.\ 1, 2, 3, 4). As in the last case,
  in all of these cases $\psi$ is defined by the half anticanonical
  system, which is impossible.
\end{enumerate}
\end{proof}

\section{Fano threefolds with canonical singularities}
\setcounter{equation}{0}

\noindent Let $Y$ be a Gorenstein Fano threefold with at most
canonical singularities and pseudo-index $\iota_Y > 1$. Then there exists
a partial crepant resolution of singularities
 \[\psi\colon X \lra Y,\]
such that $X$ is a Gorenstein threefold with at most terminal
and $\KQ$-factorial singularities and $K_X = \psi^*K_Y$ (cf.\ 
\cite{KoMo}, Theorems 6.23 and 6.25). Hence $X$ is almost Fano and 
 \[\iota_X = \iota_Y > 1.\] 

\begin{proof}[Proof of Theorem~\ref{almFano}, Fano case.]
Let $Y$ be a Gorenstein Fano threefold with at
most canonical 
singularities with $\iota_Y > 1$, and
let $X \to Y$ be a partial crepant
resolution as described above. Applying Proposition
\ref{structure} to $X$ we get a chain
\begin{equation} \label{chain3}
  \xymatrix{X_0 \ar[r]^{\phi_0} & X_1 \ar[r]^{\phi_1} &
   \cdots \ar[r]^{\phi_{m-1}}  &  X_m \ar[r]^{\phi} & W}
 \end{equation}
where $X_0$ is another 
partial crepant resolution of $Y$, each $\phi_i$ is a
blowup of a smooth point $p_{i+1}\in X_{i+1}$,
$\phi$ is an elementary extremal 
contraction of fiber type, and $\rho({X_m})=\rho(W)+1\le 3$.

For any $i=0,\dotsc,m$ denote by 
\[\psi_i\colon X_i \lra Y_i\] 
an anticanonical model of $X_i$ (in particular $Y_0=Y$). 

Now fix $i\in\{0,\dotsc,m-1\}$ and consider $\phi_i\colon X_{i}\to
  X_{i+1}$. 
Since $-K_{X_i}$ is nef,
$p_{i+1}$ is not
contained in the exceptional locus of $\psi_{i+1}$, so
  $\psi_{i+1}(p_{i+1})$ 
is a
smooth point of $Y_{i+1}$. Denote the image point
of $p_{i+1}$ in $Y_{i+1}$ by $p_{i+1}$ as well and let $\hat{X}_{i} =
\Bl_{p_{i+1}}(Y_{i+1})$.

A simple computation shows that we have an induced map
$\hat{\psi}_{i}\colon X_{i} \to \hat{X}_i$, which is crepant. We arrive at the following commutative
diagram:

 \[\xymatrix{ {X_{i}} \ar[rr]^{\phi_{i}=\text{\Bl}_{p_{i+1}}}
\ar[dr]^{\hat{\psi}_{i}} \ar[dd]_{\psi_{i}} 
& & X_{i+1} \ar[dd]^{\psi_{i+1}}\\
    &  \hat{X}_{i} \ar[rd]_{\Bl_{p_{i+1}}} \ar[dl]_{\sigma_{i}} &
    \\
   Y_{i}& & Y_{i+1} }\]

Notice that $\rho(Y_{i})\leq\rho(\hat{X}_{i})=\rho(Y_{i+1})+1$, so
either $\rho(Y_{i})\leq\rho(Y_{i+1})$, or
$\rho(Y_{i})=\rho(Y_{i+1})+1$ and $\sigma_{i}$ is
an isomorphism. 

Observe also that $\sigma_{i}$ is an isomorphism if and only if for
every curve $C\subset X_{i}$ of
anticanonical degree zero, the image $\phi_{i}(C)$ still has 
anticanonical degree zero in $X_{i+1}$.

Repeating the construction above
for all  $i\in\{0,\dotsc,m-1\}$,
in the end we get a zigzag chain of crepant maps
and blowups:

 \begin{equation} \label{zigzag}
  \xymatrix{ & \hat{X}_0 \ar[dl]_{\sigma_0} \ar[dr]^{\Bl} & & \hat{X}_1
  \ar[dl]_{\sigma_1} \ar[dr]^{\Bl} & & \dots & & \hat{X}_{m-1}
  \ar[dl]_{\sigma_{m-1}} \ar[dr]^{\Bl} &\\
       {Y=Y_0} & & Y_1 & & & \dots & & & Y_m}
 \end{equation}
where $Y_m$ is an anticanonical model of $X_m$, so
$\rho(Y_m)\leq\rho(X_m)\leq 3$.

\begin{lemma} \label{sigma}
Let $X_0$ be as in \eqref{chain3}. Assume that $m\geq 1$, and that
$X_1\not\simeq\PN_3$ if $m=1$.
Then there exists a curve $C\subset X_0$
such that $-K_{X_0}. C=0$ and $-K_{X_1}.\phi_0(C)>0$.
\end{lemma}
Using this Lemma, we complete the proof. Indeed, 
if $m=0$, then $\rho(Y)\le\rho(X_0)\leq 3$.
If $m>0$,
applying Lemma \ref{sigma} to $X_0,\dotsc,X_{m-1}$ 
in \eqref{chain3}, we see
that:
\begin{enumerate}
\item for every $i<m-1$, $\sigma_{i}$ is not an isomorphism, so
  $\rho(Y_i)\leq\rho(Y_{i+1})$;
\item for $i=m-1$, either $\sigma_{m-1}$ is not an isomorphism and
$\rho(Y_{m-1})\leq\rho(Y_m)\leq 3$, or $Y_m=X_m\simeq \PN_3$ and
$Y_{m-1}=X_{m-1}\simeq \Bl_{p_m}(\PN_3)$, so $\rho(Y_{m-1})=2$.
\end{enumerate}
This shows
$$ \rho(Y)\leq\rho(Y_{m-1})\leq 3.$$

\vspace{0.2cm}

Assume now that  $\rho(Y)=3$. Then $\rho(Y_{m-1})\geq 3$, so by 2.)
above we get
$$ 3\leq \rho(Y_{m-1})\leq \rho(Y_{m})\leq \rho(X_{m})=\rho(W)-1\leq 3.$$
We must have everywhere equality, i.e., $\rho(W) = 2$ and
$\rho(Y_m) = \rho(X_m) = 3$.  This implies that $Y_m$ is a Fano
$\PN_1$-bundle over a smooth surface, hence it is smooth with $\iota_{Y_m}=\iota_{X_m}>1$,  and we get
$Y_m=\PN_1 \times \PN_1 \times \PN_1$.

Suppose that $m>0$. Then
$X_{m-1} = \Bl_{p_m}(\PN_1 \times \PN_1 \times \PN_1)$, and
its anticanonical model $Y_{m-1}$ has Picard number
one. This contradicts $\rho(Y_{m-1})=3$, so $m=0$ and $Y=\PN_1 \times
\PN_1 \times \PN_1$. 
\end{proof}

\begin{proof}[Proof of Lemma~\ref{sigma}.]
We keep the same notation as in the proof of Theorem \ref{almFano}.

Suppose that 
there is $j\geq 1$ such that
$p_1$ lies on the strict transform $\hat{E}_1\subset X_1$  
of the exceptional
divisor $E\subset X_j$ of $\phi_j$.

For every $i=1,\dotsc,j$ denote by $\hat{E}_i$ the strict transform of
$E$ in $X_i$. We claim that $p_i\not\in\hat{E}_i$ for every
$i=2,\dotsc,j$, so that $\hat{E}_1\simeq E\simeq\PN_2$.
Indeed by the adjunction formula, $E$ is covered by
rational curves of anticanonical degree $2$. 
If we blow up a further
point $p_i\in \hat{E}_i$, then 
$\hat{E}_{i-1}\simeq \PF_1$, and
 the fibers of $\PF_1$ are now anticanonically trivial in $X_{i-1}$. 
So if we blow up any other point on $\hat{E}_{i-1}$, 
the anticanonical bundle will
not be nef any more. This implies $i=1$.
 
Then any line in $\hat{E}_1$
through $p_1$ has anticanonical degree two in $X_1$, while its proper
transform in $X_0$ is anticanonically trivial.

\vspace{0.2cm}

Suppose now that $p_1$ does not lie on the strict transform of the
exceptional divisor of any $\phi_j$.

We claim that it is enough to prove the following:
\emph{for any point $q\in X_m$ there exists a curve $C\subset
X_m$ containing $q$ and such that $-K_{X_m}. C\leq 3$}. 
In fact, choose such a curve $C$ through the image of
$p_1$ in $X_m$, and let $F\subset X_0$
be the exceptional divisor of $\phi_0$.
Let $\hat{C}_1$ and $\hat{C}_0$
 be the strict transforms of $C$ in $X_1$ and $X_0$ respectively. Then
$$-K_{X_0}.\hat{C}_0=-K_{X_1}. \hat{C}_1-2 F.
 \hat{C}_0\leq -K_{X_1}. \hat{C}_1-2\leq -K_{X_m}. {C}-2\leq 1,$$ 
so $\hat{C}_0$ is anticanonically trivial while
$-K_{X_1}. \hat{C}_1\geq 2$.

\vspace{0.2cm}

Consider now the contraction $\phi\colon X_m\to W$.

\vspace{0.2cm}

1.) Assume that $\dim W = 0$, i.e. $X_m$ is a Fano threefold with $\rho(X_m) =
1$.   

\smallskip

If $\iota_{X_m} \ge 4$, then by Proposition \ref{description} we have
$X_m \simeq \PN_3$, so by hypothesis $m\geq 2$.
 Observe that $X_{m-1}\simeq \Bl_{p_m}(\PN_3)$ 
is covered by curves of anticanonical degree two, so we
are done (just replace $X_m$ by $X_{m-1}$, this is
possible because $m\geq 2$).

\smallskip

If $\iota_{X_m} = 3$, Proposition \ref{description} says that
$X_m$ is a quadric; in particular it is covered by curves of
anticanonical degree $3$.

\smallskip

Assume that $\iota_{X_m}=2$. We show that through any point of $X_m$
there is a curve of anticanonical degree two. Since $\rho(X_m) = 1$,
we have $r_{X_m} = \iota_{X_m}$ by Lemma~\ref{index}.

Let $\pi\colon \X \to \Delta$ be the smoothing of $X_m$. Then the
general fiber $\X_t$ is a smooth Fano threefold with $\rho(\X)t)=1$
and  $r_{\X_t} = \iota_{\X_t} = 2$ (see Theorem \ref{Smoothing}).
Let $x \in X_m$ be any point and $x_t \in \X_t$, $t \not= 0$ be points
with limit $x$. There exist
rational curves $\C_t$ in $\X_t$
containing $x_t$, such that $-K_{\X_t}.\C_t =
r_{\X_t} = r_{X_m}$ is constant. Let $\C_0$ be the limit curve. Then
$-K_{X_m}.\C_0 = r_{X_m} = \iota_{X_m}$ and $\C_0$ contains $x$. If
$\C_0$ is not irreducible, then $-K_{X_m}.\C_{0,i} < \iota_{X_m}$ for
some component, which is impossible. Hence $\C_0$ is an irreducible
rational curve.

\vspace{0.2cm}

2.) Assume that $\dim W = 1$, i.e., $W \simeq \PN_1$ and $X_m \to W$ is a
del Pezzo fibration with general fiber $F$. 
Since $-K_{X_m}|_F = -K_F$, we have $\iota_F > 1$,
hence $F \simeq 
\PN_2$ or $\PN_1\times \PN_1$. 
In both cases $X_m$ is covered by  rational curves
contained in the fibers of $\phi$, having
anticanonical degree at most three.

\vspace{0.2cm}

3.) Finally, when
$\dim W = 2$, $X_m$ is a $\PN_1$--bundle and it
is again covered by rational curves of anticanonical degree two. This
finishes the proof.
\end{proof}

\begin{proof}[Proof of Theorem \ref{Mukcan}]
The statement 
is trivial for $\iota_Y = 1$, it follows from Theorem~\ref{almFano} for $\iota_Y
= 2$, and from Proposition
\ref{description}  for $\iota_Y \ge 3$.
\end{proof}

\section{The toric case}
\setcounter{equation}{0}

\noindent In order to prove Proposition \ref{toric}, we first need the following bound.
\begin{lemma}
\label{24}
Let $X$ be a smooth toric almost Fano threefold of index $r_X=2$. If
the anticanonical model of
$X$ is not $\KQ$-factorial,
then $(-K_X)^3\geq 24$.
\end{lemma}
\begin{proof}
There exists $L\in\Pic(X)$ such that $-K_X=2L$.
Recall that on a smooth toric variety, every nef line bundle is
globally generated (see for instance \cite{torimori}, Lemma on p.\
261). Hence
$L$ is globally
generated
and big, and it defines a map $\varphi_{L}\colon X\to
\PN_N$, where $N:=h^0(L)-1$. Since $\varphi_L$ contracts all
anticanonically trivial curves, it 
factors through the anticanonical model $Y$ of $X$:
$$\xymatrix{
X\ar[rr]^{\varphi_{L}}\ar[dr]_{\psi} && {\PN_N}\\
& Y \ar[ur]&
}$$
Notice that $-K_Y$ is the pull-back of $\mathcal{O}_{\PN_N}(2)$, so
$r_Y=2$. 

Set $Z:=\varphi_{L}(X)\subset\PN_N$, and observe that $L^3=(\deg
Z)(\deg\varphi_L$).

Since $(-K_X)^3=8L^3$, we have to show that $L^3\geq 3$.
By contradiction, if $L^3< 3$, we get the following possibilities:
\begin{enumerate}
\item
$L^3=1$, $Y\simeq Z=\PN_{3}$: this is impossible because $Y$ has index two;
\item
$L^3=2$, $Y\simeq Z$ a quadric in $\PN_{4}$: again, this is impossible
because $Y$ has index two;
\item
$L^3=2$, $Z=\PN_{3}$, $Y\to \PN_{3}$ an equivariant finite map. This
means that the fan of $Y$ is the same as the fan of $\PN_{3}$, with
respect to a sublattice of $\mathbb{Z}^3$. Hence every cone of the fan
is simplicial and $Y$ is $\KQ$-factorial, a contradiction. 
\end{enumerate}
\end{proof}
\begin{proof}[Proof of Proposition \ref{toric}]
Recall that an almost Fano toric threefold always admits a crepant
toric resolution $X'\to X$ (see \cite{nill}, Proposition 1.15). 
Hence $X'$ is a smooth almost Fano with
$\rho(X')\ge \rho(X)$, and it is enough to prove the statement in the
smooth case.

So assume that $X$ is smooth. Applying Proposition \ref{structure}, we get a
diagram as \eqref{fact},
so that 
$$\rho(X)=\rho(X_0)=\rho(X_m)+m, \quad\text{and}\quad
0<(-K_X)^3=(-K_{X_m})^3-8m.$$ 
Recall that $X_0$ is obtained from $X_m$ by $m$ blowups. 
 Since $X$ is smooth and
toric, the same holds for all varieties in \eqref{fact}, and the
maps are equivariant. In particular, the center of each blowup must be
a fixed point for the torus action.
Recall also that any
elementary contraction of fiber type between smooth toric varieties is
a projective bundle. 

If $W=\PN_1$, then $X_m$ is a $\PN_2$-bundle over $\PN_1$.
Through any point of $X_m$ there is a
curve of anticanonical degree 3, so any blowup of $X_m$ will contain a
curve of anticanonical degree 1. This implies $m=0$ and $\rho(X)=2$.

Suppose that $\dim W\neq 1$.
If $m=0,1$ we have $\rho(X)=\rho(X_m)+m\leq 4$,
 so we can assume that
$m\geq 2$.
Then Proposition \ref{structure} yields that $W$ is either a point or
$\PN_2$. 

This implies that  $X_m$ is either $\PN_3$ or
$\PN_{\PN_2}(\mathcal{O}\oplus\mathcal{O}(3))$. 
Indeed
 $X_m\simeq\PN_3$ when $W$ is a point. If $W\simeq\PN_2$, then
$X_m$
is a $\PN_1$-bundle over $\PN_2$.  
Excluding 
$\PN_1\times\PN_2$, which 
has been considered above, and $\Bl_p \PN_{3}$, 
the only possibility is $\PN_{\PN_2}(\mathcal{O}\oplus\mathcal{O}(3))$.

In both cases, $-K_{X_m}$ is divisible by two in $\Pic(X_m)$. This implies that
$-K_{X_0}$ is divisible by two in $\Pic(X_0)$, hence $2|r_{X_0}$. On the
other hand, $\rho(X_0)>1$, so $r_{X_0}=2$.

Observe that
$X_{m-2}$ contains at least 
one curve
$C\simeq\PN_1$ with normal bundle
$\mathcal{O}_{\PN_1}(-1)^{\oplus 2}$, such that $C$ is not contained
in a surface covered by anticanonically trivial curves (such surface
should be isomorphic to $\mathbb{F}_1$).
Hence the
the same holds for $X_0$, and this implies that the anticanonical model
of $X_0$ is not $\KQ$-factorial.

Now Lemma \ref{24} yields
$(-K_X)^3=(-K_{X_0})^3\geq 24$, hence
$$m\leq\frac{1}{8}(-K_{X_m})^3-3.$$

Therefore if $W$ is a point and $X_m\simeq\PN_3$, we have $m\leq 5$
and $\rho(X)\leq 6$.
If $W\simeq\PN_2$ and
$X_m=\PN_{\PN_2}(\mathcal{O}\oplus\mathcal{O}(3))$, we have
$(-K_{X_m})^3=72$, so $m\leq 6$ and
$\rho(X)\leq 8$. 

\medskip

Finally, one can check directly that there are choices of 6 blowups
of fixed points on $\PN_{\PN_2}(\mathcal{O}\oplus\mathcal{O}(3))$ such
that the resulting variety is an almost Fano threefold with Picard
number 8.
\end{proof}


\end{document}